\documentclass[12pt]{article}

\usepackage[english]{babel}
\usepackage[utf8x]{inputenc}
\usepackage[T1]{fontenc}

\usepackage[letterpaper,top=1in,bottom=1in,left=1in,right=1in,marginparwidth=1.75cm]{geometry}

\usepackage{amsmath,amsfonts,amsthm,amssymb,mathrsfs,dsfont} 
\usepackage{mathtools}
\usepackage{graphicx}
\usepackage[colorinlistoftodos]{todonotes}
\usepackage[colorlinks=true, allcolors=blue]{hyperref}
\usepackage[capitalize,nameinlink]{cleveref}
\usepackage{verbatim}
\usepackage{enumerate}
\usepackage{slantsc}

\newcommand{\R}{\mathbb{R}}

\global\long\def\R{\mathbb{R}}

\global\long\def\S{\mathbb{S}}

\global\long\def\P{\mathbb{P}}

\global\long\def\F{\mathcal{F}}

\global\long\def\E{\mathbb{E}}

\global\long\def\S{\mathcal{S}}

\newtheorem{theorem}{Theorem}[section]
\newtheorem*{namedtheorem}{\theoremname}
\newcommand{\theoremname}{testing}

\newtheorem{thm}[theorem]{Theorem}

\newtheorem{lem}[theorem]{Lemma}
\newtheorem{claim}[theorem]{Claim}

\newtheorem{cor}[theorem]{Corollary}

\newtheorem{conjecture}[theorem]{Conjecture}
\newtheorem*{question*}{Question}

\theoremstyle{definition}
\newtheorem{definition}[theorem]{Definition}

\newtheorem{remark}[theorem]{Remark}

\theoremstyle{plain}

\title{Nondegenerate spheres in four dimensions}
\author{Thao~Do \footnote{Massachusetts Institute of Technology, Department of Mathematics, thaodo@mit.edu} }

\begin{document}
\maketitle

\begin{abstract}
Non-degeneracy was first defined for hyperplanes by Elekes-T\'oth in \cite{elekes2005incidences}, and later extended to spheres by Apfelbaum-Sharir in \cite{apfelbaum2011non}:
given a set $P$ of $m$ points in $\R^d$ and some $\beta\in(0,1)$, a $(d-1)$-dimensional sphere (or a $(d-1)$-sphere) $S$ in $\R^d$ is called \emph{$\beta$-nondegenerate} with respect to $P$ if $S$ does not contain a proper subsphere   $S'$ such that $|S'\cap P|\geq \beta|S\cap P|$. 
Apfelbaum-Sharir found an upper bound for the number of incidences between points and nondegenerate spheres in $\R^3$, which was recently used by Zahl in \cite{zahl2017breaking} to obtain the best known bound for the unit distance problem in three dimensions. 

In this paper, we show that the number of  incidences between $m$ points and $n$ $\beta$-nondegenerate 3-spheres in $\R^4$ is $O_{\beta,\varepsilon}\left(m^{\frac{15}{19}+\varepsilon} n^{\frac{16}{19}}+mn^{\frac{2}{3}}\right)$.
As a consequence, we obtain a bound of $O_{\varepsilon}(n^{2+4/11+\varepsilon})$ on the number of similar triangles formed by $n$ points in $\R^4$, an improvement over the previously best known bound $O(n^{2+2/5})$. 

While proving this, we find it convenient to work with a more general definition of nondegeneracy: a bipartite graph $G=(P,Q)$ is called \emph{$\beta$-nondegenerate} if $|N(q_1)\cap N(q_2)|<\beta |N(q_1)|$ for any two distinct vertices $q_1,q_2\in Q$; here $N(q)$ denotes the set of neighbors of $q$  and $\beta$ is some positive constant less than 1.  A $\beta$-nondegenerate graph can have up to $\Theta(|P||Q|)$ edges  without any restriction, but must have much fewer edges if the graph is \emph{semi-algebraic} or has bounded \emph{VC-dimension}. 
We show that results in \cite{elekes2005incidences} and \cite{apfelbaum2011non} still hold under this new definition, and so does our new bound for spheres in $\R^4$.
\end{abstract}

\section{Introduction}
\subsection{Nondegenerate hyperplanes and spheres}
Counting the number of incidences between points and algebraic varieties (such as hyperplanes or spheres) in high dimensions has been a center of interest in incidence geometry. However, this problem is only interesting under certain constraints -- indeed, an arbitrary set of $m$ points and a set of $n$ hyperplanes in $\R^d$ for $d\geq 3$ can have up to $mn$ incidences by putting all the points in a low dimensional flat and letting all the hyperplanes contain that flat. To avoid this situation, one can impose conditions such as the points must be vertices in the arrangement of the hyperplanes (see \cite{agarwal1992counting}), or the incidence graph is $K_{r,s}$-free for some fixed small $r,s$ (Zarankiewicz's problem). In this paper, we focus on another type of restriction: nondegeneracy.

Nondegeneracy was introduced by Elekes and T\'oth 
in \cite{elekes2005incidences}. 
Given a set $P$ of $m$ points in $\R^d$, a hyperplane $H$ is \emph{$\beta$-nondegenerate} if there does not exist a lower dimensional flat $F\subset H$ that contains more than $\beta$ fraction of the number of points of $P$ in $H$, i.e. $|F\cap P|>\beta|H\cap P|$. Otherwise, $H$ is $\beta$-degenerate. Elekes and T\'oth proved \footnote{What they actually proved is that the maximum number of $\beta$-nondegenerate, $k$-rich (i.e. containing at least $k$ points of $P$) hyperplanes is  
$O_{\beta,d} \left(\frac{m^{d+1}}{k^{d+2}}+\frac{m^{d-1}}{k^{d-1}}\right)$
for $\beta<\beta_d$ for some small $\beta_d$. It is later shown in \cite{lund2017two} that we can indeed let $\beta_d=1$. See \cite{apfelbaum2007large} for how this implies \eqref{eq:E-T}.} 
that $m$ points and $n$ $\beta$-nondegenerate hyperplanes in $\R^d$ can have at most
\begin{equation}\label{eq:E-T}
O_{d,\beta}\left((mn)^{\frac{d}{d+1}}+mn^{1-\frac{1}{d-1}}\right)
\end{equation}
incidences. Here we use the notation $O_{d}$ to indicate the constant depends on $d$.

The concept of nondegeneracy extends naturally to spheres. Given a set $P$ of $m$ points in $\R^d$ and some $\beta\in(0,1)$, a $(d-1)$-dimensional sphere (or a $(d-1)$-sphere) $S$ in $\R^d$ is \emph{$\beta$-nondegenerate} with respect to $P$ if it does not contain a proper subsphere   $S'$ such that $|S'\cap P|\geq \beta|S\cap P|$. 
By lifting the spheres to be hyperplanes\footnote{via the map $(x_1,\dots, x_d)\mapsto (x_1,\dots, x_d, x_1^2+\dots+x_d^2)$} in $\R^{d+1}$, we immediately get the following bound on the number of incidences between $m$ points and $n$ $\beta$-nondegenerate spheres in $\R^d$:
\begin{equation}\label{eq:lifting sphere to hp}
O_{d,\beta}\left((mn)^{\frac{d+1}{d+2}}+mn^{1-\frac{1}{d}}\right).
\end{equation}

When $d=3$, Apfelbaum and Sharir in \cite{apfelbaum2011non}  improved this bound to\footnote{To see this bound is indeed an improvement from  $(mn)^{4/5}+mn^{2/3}$: clearly $mn^{1/2}\leq mn^{2/3}$; $m^{8/11}n^{9/11}\leq (mn)^{4/5}$ iff $n\leq m^4$, which holds true by \cref{lem:spanning}.}
\begin{equation}\label{eq:Ap-Sh}
O_\varepsilon\left(m^{\frac{8}{11}+\varepsilon}n^{\frac{9}{11}}+mn^\frac{1}{2}\right).
\end{equation}

These incidence bounds with nondegenerate hyperplanes and spheres have many applications in discrete geometry. For example, the bound \eqref{eq:E-T} was used in \cite{afshani2016applications} to derive an efficient algorithm for the point covering problem in $\R^3$. It was also used in \cite{apfelbaum2007large} to derive a lower bound on the maximum size of a complete bipartite subgraph one can find in a point-hyperplane incidence graph with many edges. Nondegenerate spheres have been used in the analysis of the distinct distance problem in \cite{aronov2003distinct} and in the mutually similar simplices problem in \cite{agarwal2007similar}. Recently,  Zahl \cite{zahl2017breaking} used \eqref{eq:Ap-Sh} to obtain a bound of $O(n^{295/197 + \epsilon})$ on the number of unit formed by $n$ points in $\R^3$, the first improvement over the well-known $n^{3/2}$ upper bound.

It is not obvious how to extend the concept of nondegeneracy to other geometric objects. In this paper, we propose a similar but slightly different definition of nondegeneracy that can be applied more widely. Under this more general framework, we extend Apfelbaum-Sharir's method in \cite{apfelbaum2011non} to some other families of 
surfaces in $\R^3$ and use that to prove a new upper bound for incidences with nondegenerate spheres in four dimensions.

\subsection{New definition of nondegeneracy} 
\begin{definition}\label{def:nondeg graph}
Consider a bipartite graph $G=(P,Q)$ and some $0<\beta<1$. For each $q\in Q$, let $N(q)$ denote the set of neighbors of $q$ in $P$. We say a vertex $q\in Q$ is \emph{$\beta$-nondegenerate} with respect to $P$ if $|N(q)\cap N(q')|<\beta |N(q)|$ for any $q'\in Q, q'\neq q$ (otherwise $q$ is $\beta$-degenerate). We say $G$ is a \emph{$\beta$-nondegenerate graph} if every vertex of $Q$ is $\beta$-nondegenerate.
\end{definition}
\begin{remark}\label{rmk:nondeg graph}
\begin{enumerate}[i]
\item This definition is not symmetric between $P$ and $Q$.  We sometimes make it clear by saying $Q$ is $\beta$-nondegenerate with respect to $P$. If $P$ is also $\beta$-nondegenerate with respect to $Q$, we say $G$ is \emph{dually $\beta$-nondegenerate}.  
\item\label{rmk:Zarankiewicz} If all vertices in $Q$ have $d$ neighbors (which holds for $d$-regular graphs), the condition of $\beta$-nondegeneracy is equivalent to $G$ being $K_{\beta d,2}$-free. This is related to Zarankiewicz's problem, that is, finding the maximum number of edges in a $K_{r,s}$-free graph on $n$ vertices where $r,s$ are fixed.
\item Comparison with the previous definitions of nondegenerate hyperplanes and spheres: It is easy to see that the incidence graphs between points and $\beta$-nondengerate hyperplanes (or spheres) in some Euclidean space are $\beta$-nondegenerate graphs. Indeed, since a hyperplane (sphere) does not contain any subflat (subsphere) of more than $\beta$ fraction number of points of $P$, its intersection with any other hyperplane (sphere) must contain fewer than a $\beta$ fraction of the points. 
The other direction does not always hold, hence this definition strictly generalizes the previous one. Therefore, any upper bound for points-hyperplanes or points-spheres incidence under this new definition would immediately imply the same bound for the old definition. 
\end{enumerate}

\end{remark}
From now on, we shall use this definition of nondegeneracy, and refer to the old one as \emph{geometrically nondegenerate}.

\subsection{Our results}
The questions we would like to investigate have the following form: what is the maximum number of edges a graph $G = (P,Q)$ can have, given that $G$ is $\beta$-nondegenerate and satisfies various other conditions?  
It turns out that $\beta$-nondegeneracy alone is not enough for a $o(|P||Q|)$ bound on $|E(G)|$:
\begin{thm}\label{thm:general graph}
For any $0<\beta<1$, there exists a $\beta$-nondegenerate bipartite graph $G=(P,Q)$ with at least $\frac{\beta}{6} |P||Q|$ edges.
\end{thm}
However, as we will show, stronger bounds on $|E(G)|$ hold for graphs which are \emph{semi-algebraic} and graphs with bounded \emph{VC dimension}. When $G$ is an incidence graph between points and certain families of geometric objects such as planes and spheres in $\R^3$, we show that a bound similar to \eqref{eq:E-T} and \eqref{eq:Ap-Sh} holds. Finally, our main result is a new upper bound on the number of incidences between nondegenerate 3-spheres in $\R^4$, which implies an improved bound on the number of similar triangles spanned by a collection of points in $\R^4$.

\subsubsection{Graphs with bounded VC dimensions and semi-algebraic graphs}
A \emph{set system} $F$ over a ground set $P$ is  a collection of subsets of $P$ (here we allow $F$ to contain repeated elements). The \emph{Vapnik-Chervonenkis dimension (VC-dimension)} of $F$ is the largest integer $d$ for which there exists a $d$-element set $S\subset P$ such that for every $B$ subset of $S$, we can find a set $A\in F$ such that $A\cap S=B$. The \emph{VC-dimension of a graph} with vertex set $V$ is the VC-dimension of the set system formed by the neighbor sets over the ground set $V$. When $G$ is bipartite, say $G=(P,Q)$, we can define two set systems:  $F_1$ consists of all neighbor sets of $Q$ over the ground set $P$, and $F_2$ consists of all neighbor sets of $P$ over the ground set $Q$. We use $left$ and $right$ $VC$-dimension to denote the VC-dimensions of $F_1$ and $F_2$ respectively. It is easy to see that the VC-dimension of $G$ is the maximum of its left and right VC-dimensions.

Graphs with bounded VC-dimension obey a dramatically improved regularity lemma \cite{lovasz2010regularity,alon2007efficient}: here $\epsilon$-regular partitions require only poly$(1/\epsilon)$ parts, as opposed to the tower-type bound which is tight for general graphs. Using this, Fox, Pach and Suk in \cite{fox2017erdos} get improvements on the Erd\H{o}s-Hajnal conjecture bound and Turan-Ramsey number bound for graphs with bounded VC dimensions. In this paper we show another special property of such graphs.
\begin{thm}\label{thm:VC dim}
Assume a bipartite graph $G=(P,Q)$ has left-VC-dimension at most $d$ and $Q$ is $\beta$-nondegenerate with respect to $P$, then $|E(G)|=O_{d,\beta}( mn^{1-1/d}+n)$ where $m=|P|,n=|Q|$. 
\end{thm}
This implies a similar result for \emph{semi-algebraic graphs}. 
Fixing some positive integers $d_1,d_2$, let $G=(P,Q,\mathcal{E})$ be a bipartite graph on sets $P$ and $Q$, where we think of $P$ as a set of $m$ points in $\R^{d_1}$ and $Q$ as a set of $n$ points in $\R^{d_2}$. We say $G$ is \emph{semi-algebraic with description complexity $t$} if there are $t$ polynomials $f_1,\dots, f_t\in \R[x_1,\dots, x_{d_1+d_2}]$, each of degree at most $t$ and a Boolean function $\Phi(X_1,\dots, X_t)$ such that for any $p\in P, q\in Q$: $$ (p,q)\in \mathcal{E} \iff \Phi(f_1(p,q)\geq 0,\dots, f_t(p,q)\geq 0)=1.$$
In other words, we can describe the incidence relation using at most $t$ inequalities involving polynomials of degree at most $t$. Semi-algebraic graphs have been  studied recently in \cite{fox2014density,fox2016polynomial,conlon2014ramsey}. They are important in discrete geometry because many graphs coming from discrete geometry problems (such as incidence graphs and intersection graphs of popular geometric objects) are semi-algebraic. As shown in \cite{fox2014semi},  such a graph satisfies a same condition on the \emph{shatter function}\footnote{defined at the beginning of the proof of \cref{thm:VC dim}.} that a graph with VC dimension $d_2$ satisfies. Therefore, we have the following result. 
\begin{thm}\label{thm:semi-alg graph}
Assume $G=(P,Q)$ is a semi-algebraic graph with description complexity $t$ as above. If $G$ is $\beta$-nondegenerate, then it has at most $O_{\beta, d_1,d_2, t}(mn^{1-1/d_2}+n)$ edges.
\end{thm}
This immediately implies if $G$ is the incidence graph between $m$ points and $n$ $\beta$-nondegenerate hyperplanes in $\R^d$, then $|e(G)|=O(mn^{1-1/d}+n)$. Similarly, if $G$ is the incidence graph between $m$ points and $n$ $\beta$-nondegenerate $(d-1)$-spheres in $\R^d$, then $|e(G)|=O(mn^{1-1/(d+1)}+n)$. However, as mentioned earlier, we expect much stronger bounds for these cases.
\subsubsection{Nondegenerate surfaces in $\R^3$}
We show that the Apfelbaum-Sharir bound \eqref{eq:Ap-Sh} can be extended to other families of surfaces in $\R^3$ that satisfy certain constraints. 
\begin{definition}\label{def:good family}
Let $P$ be a set of finite points in $\R^3$. Given positive integers $s,t$, we say a $\F$ is a \emph{good $s$-dimensional family with respect to $P$} with complexity $t$ in $\R^3$ if 
\begin{enumerate}[i]
\item Each element of $\F$ is an algebraic surface in $\R^3$ of degree at most $t$.
\item The parameter space of $\F$ is a $s$-dimensional subvariety of $\R^{t+3\choose 3}$, the space of all degree $t$ polynomials in three variables.
\item Any two surfaces in $\F$ intersect at at most one irreducible curve.
\item Let $\Gamma$ denote the set of all intersection curves between two surfaces in $\F$, then $\Gamma$ belongs to a family of curves with $s-1$ degrees of freedom and multiplicities at most $t$, i.e. for any tuple of $s-1$ points in $P$ there are at most $t$ curves in $\Gamma$ passing through all of them. 
\end{enumerate}
\end{definition} 
 For example, planes in $\R^3$ form a good 3-dimensional family and spheres in $\R^3$ form a good 4-dimensional family. Given a set of points $P$ in $\R^3$, a surface $S$ in $\F$ is called geometrically $\beta$-nondegenerate (w.r.t. P) if there does not exist a curve $\gamma\in \Gamma$ such that $\gamma\subset S$ and $|\gamma\cap P|\geq\beta|S\cap P|$. 
 
\begin{thm}\label{thm:extend Ap-Sh in R3}
The number of incidences between a set $P$ of $m$ points and a set $Q$ of $n$ $\beta$-nondegenerate algebraic surfaces taken from a good $s$-dimensional family $\F$ in $\R^3$ is $O_{s,t,\beta,\varepsilon}\left(m^{\frac{2s}{3s-1}+\varepsilon}n^{\frac{3s-3}{3s-1}}+mn^{1/2}+n\right).$
\end{thm}

Note that the term $n$ is unnecessary here if $Q$ is a set of geometrically nondegenerate planes or spheres because of the following claim (for the proof see section \ref{sec:R3}).

\begin{claim}\label{lem:spanning}
Let $P$ be a set of $m$ points in $\R^d$ and $\beta$ some constant between 0 and 1. There are at most ${m\choose d}$ geometrically $\beta$-nondegenerate hyperplanes with respect to $P$. Similarly, there are at most ${m\choose d+1}$ geometrically $\beta$-nondegenerate spheres with respect to $P$.
\end{claim}

Indeed, by the claim, for planes in $\R^3$, $n\leq m^3$, which implies $n\leq (mn)^{3/4}$. For spheres, $n\leq m^4$ which implies $n\leq m^{8/11}n^{9/11}$.
This result is thus a generalization of  \eqref{eq:Ap-Sh} and \eqref{eq:E-T} (in three dimensions and up to $\varepsilon$).
It is not easy to find any good family of surfaces other than planes and spheres, hence this theorem may appear pointless. However, we include it here because its proof implies a result which we will need in the proof of our main result, \cref{thm:nondeg spheres 4d}. 

\begin{cor}\label{cor:projected sphere}
Let $P$ be a set of $m$ points and $Q$ a set of $n$ $\beta$-nondegenerate 2-spheres in $\R^d$ for $d\geq 3$. Let $I(P,Q)$ denote the number of incidences between $P$ and $Q$, then 
\begin{equation}\label{eq:projected sphere}
O_{d,\varepsilon}\left(m^{\frac{8}{11}+\varepsilon}n^{\frac{9}{11}}+mn^\frac{1}{2}+n\right).
\end{equation}
\end{cor}
Basically this means \eqref{eq:Ap-Sh} still holds for nondegenerate 2-spheres in any dimension $d$ at least 3.
We prove this by projecting to a generic 3-dimensional subspace, and think of the projected 2-spheres as a good 4-dimensional family.

One can find the proofs of these three results in section \ref{sec:R3}. Before the proof of \cref{thm:extend Ap-Sh in R3}, we sketch the main idea and how it is related to the proof for spheres in \cite{apfelbaum2011non}.
\subsubsection{Nondegenerate spheres in $\R^4$}
This is our main new result. 

\begin{thm}\label{thm:nondeg spheres 4d}
The number of incidences between $m$ points and $n$ (geometrically) $\beta$-nondegenerate $3$-spheres in $\R^4$ is $O_{d,\beta,\varepsilon}\left(m^{\frac{15}{19}+\varepsilon} n^{\frac{16}{19}}+mn^{\frac{2}{3}}+n\right).$
\end{thm}

\begin{remark}
\begin{enumerate}[i]
\item If assuming geometrically nondegenerate, we can ignore the last term.
\item This  bound is stronger than $O(m^{5/6}n^{5/6}+mn^{3/4})$, the  bound in \eqref{eq:lifting sphere to hp} for $d=4$.
\item By a result in \cite{fox2014semi}, the number of incidences between $m$ points and $n$ 3-spheres in $\R^4$ without $K_{r,s}$ in their incidence graph is $O_\varepsilon(m^{\frac{15}{19}+\varepsilon} n^{\frac{16}{19}}+m+n)$. Notice the first term here is the same with the first term in \eqref{eq:sphere 4d}. 
See also \cref{rmk:nondeg graph}\ref{rmk:Zarankiewicz}.

\item The second term $mn^{2/3}$ is the same with the second term in Elekes-T\'oth's bound \eqref{eq:E-T} for hyperplanes in $\R^4$. This is consistent with Apfelbaum-Sharir's result \eqref{eq:Ap-Sh}. 

\end{enumerate}
\end{remark}
This theorem is proved in section \ref{sec:R4}. Before the proof, we give a brief outline and how it is related to the proof of \cref{thm:extend Ap-Sh in R3}.

We conjecture that a similar result holds for spheres in all dimensions. In the bound below, the first term is the same with the one in \cite{fox2014semi} for point-sphere incidences in $\R^d$ with no $K_{r,s}$ for big $r,s$, and the second term is the same with the second one in \eqref{eq:E-T}.
\begin{conjecture}
Given $d\geq 3$ and $\beta\in(0,1)$. The number of incidences between $m$ points and $n$ (geometrically) $\beta$-nondegenerate $(d-1)$-spheres in $\R^d$ is $O_{\varepsilon,d}(m^{\frac{d^2-1}{d^2+d-1}+\varepsilon} n^{\frac{d^2}{d^2+d-1}}+mn^{1-\frac{1}{d-1}}+n).$
\end{conjecture}
\subsubsection{Application to the similar triangle problem}
Erd\H{o}s and Purdy asked the following question (see Problem 10 section 6.1 in \cite{brass2006research}): What is the maximum number of mutually similar $k$-dimensional simplices among $n$ points in $\R^d$? Let denote that number by $f_{k,d}(n)$. When $k=2$, it is called the similar triangles problem. As mentioned in \cite{agarwal2007similar},  the only interesting cases are $d=3,4$ and 5 because  $f_{2,2}(n)=\Theta(n^2)$ and $f_{2,6}=\Theta(n^3)$. In
that paper, Agarwal, Apfelbaum, Purdy and Sharir showed that 

$$f_{2,3}(n)=O(n^{13/6});\quad f_{2,4}(n)=O(n^{12/5});\quad f_{2,5}(n)=O(n^{8/3}).$$

Later in \cite{apfelbaum2011non}, Apfelbaum and Sharir used \eqref{eq:Ap-Sh} to show $f_{2,3}(n)=O_\varepsilon(n^{58/27+\varepsilon})$. In \cref{section:application}, we prove a better  bound for $f_{2,4}(n)$. 
\begin{thm}\label{thm:application to similar triangles in R4}
Let $P$ be a set of $n$ points in $\R^4$ and $\Delta$ is a fixed triangle. The number of triples of points in $P$ that span a triangle similar to $\Delta$ is at most $O_\varepsilon(n^{2+4/11+\varepsilon})$.
\end{thm}

\subsection{Organization} In the subsequent sections, we present the proofs for Theorems \ref{thm:general graph}, \ref{thm:VC dim}, \ref{thm:semi-alg graph},
\ref{thm:extend Ap-Sh in R3}, \ref{thm:nondeg spheres 4d}, and \ref{thm:application to similar triangles in R4} respectively.

\subsection{Acknowledgements} The author would like to thank Larry Guth and Micha Sharir for helpful conversations. She also thanks Ethan Jaffe, Vishesh Jain and Jake Wellens for some helpful comments.

\section{General nondegenerate graphs}
\begin{proof}[Proof of \cref{thm:general graph}]
We use the probabilistic method. Assume $|P|=m,|Q|=n$. Pick each edge in $P\times Q$ independently with probability $\rho$ for some $\rho\in (0,1)$ chosen later. For each $q,q'\in Q$, let $X_q$ be the random variable that represents the number of neighbors of $q$ and $X_{q,q'}$ represents the number of common neighbors of $q$ and $q'$ in this random graph. Each vertex $p\in P$ is a neighbor of $q$ with probability $\rho$ and is a common neighbor of $q$ and $q'$ with probability $\rho^2$. By the linearity of expectation, we have $\E[X_q]=\rho m$ and $\E[X_{q,q'}]=\rho^2 m$. Moreover, by Chernoff's bound (see for example Appendix A in \cite{alon2004probabilistic}), each variable is strongly concentrated around its expectation: 
$$\P\left(|X_q-\rho m|>\frac{\rho}{2}m\right)<e^{-\frac{1}{2}\rho^2 m}.$$
$$\P\left(|X_{q,q'}-\rho^2 m|>\frac{\rho^2}{2}m\right)<e^{-\frac{1}{2}\rho^4 m}.$$
By the union bound, the probability there exists a vertex in $Q$ of degree less than $\frac{\rho m}{2}$ or two vertices in $Q$ with more than $\frac{3\rho^2 m}{2}$ common neighbors is at most
$$me^{-\frac{1}{2}\rho^2 m}+{m\choose 2}e^{-\frac{1}{2}\rho^4 m}$$
which is less than 1 for $m$ big enough. Hence there exists a graph where $|N(q)|\geq \frac{\rho m}{2}$ and $|N(q)\cap N(q')|\leq \frac{3\rho^2 m}{2}$ for any distinct $q,q'\in Q$. This implies the graph is $\beta$-nondegnerate if we choose $\rho=\frac{1}{3}\beta$. Moreover, the graph has at least $\frac{\rho m}{2}n=\frac{\beta}{6}mn$ edges.
\end{proof}
\section{Graphs with bounded VC dimension}
\begin{proof}[Proof of \cref{thm:VC dim}]
Recall that $F_1$ is the set system $\{N(p): p\in P\}$ on the ground set $Q$ and $F_2$ is the set system $\{N(q): q\in Q\}$ on the ground set $P$. The \emph{shatter function} of $F_1$, denoted by $\pi_{F_1}$, is defined as follows: 

$$\pi(z):= \max_{B\subset Q, |B|=z}|\{A\cap B: A\in F_1\}|.$$ 

By a famous result (see for example \cite{sauer1972density,shelah1972combinatorial,vapnik2015uniform}), when $F_1$ has VC-dimension at most $d$:  $$\pi_{F_1}(z)\leq \sum_{i=0}^d {z\choose i}< cz^d $$
for some constant $c$ depending on $d$. By observation 2.6 in \cite{fox2014semi}, there exist $q_1,q_2\in Q$ such that the number of semi-algebraic sets that cross $\{q_1,q_2\}$ is at most $c'mn^{-1/d}$ for some $c'$ that depends on $c$ and $d$. This is equivalent to $|N(q_1)\Delta N(q_2)|\leq c'mn^{-1/d}$. Since $Q$ is $\beta$-nondegenerate, $|N(q_1)|\leq \frac{1}{1-\beta}|N(q_1)\setminus N(q_2)|\leq c'mn^{-1/d}$. Deleting this vertex and repeating the process, we find $$|E(G)|\leq \frac{c'}{1-\beta} \left(mn^{1-1/d}+n\right).$$
\end{proof}
\begin{proof}[Proof of \cref{thm:semi-alg graph}] As argued in \cite{fox2014semi}, the semi-algebraic graph $G$ satisfies $\pi(z)=O(z^{d_2})$ where $\pi$ is the shatter function of the set system $\{N(p):p\in P\}$. We then proceed similarly as the above proof after replacing $d$ by $d_2$.
\end{proof}
\section{Nondegenerate surfaces in $\R^3$}\label{sec:R3}
In this section we first prove \cref{thm:extend Ap-Sh in R3}, then show how its proof implies \cref{cor:projected sphere}. Finally we prove \cref{lem:spanning}.

To prove \cref{thm:extend Ap-Sh in R3}, we follow the strategy  of Apfelbaum and Sharir in \cite{apfelbaum2011non}. To extend their result for spheres to our result for surfaces of bounded degree we need  stronger tools. In particular, we need  a  \emph{cutting method} for general surfaces by Chazelle, Edelsbrunner, Guibas, and Sharir \cite{chazelle1991singly} (instead of the cutting for spheres in \cite{clarkson1990combinatorial}), incidence bound for curves with bounded degrees of freedom by Pach and Sharir \cite{pach1998number} (instead of point-circle incidence bound in \cite{marcus2006intersection}) and a bound from \cref{thm:semi-alg graph} (instead of \eqref{eq:lifting sphere to hp}).
 \\
 \\
\noindent\textit{Proof of \cref{thm:extend Ap-Sh in R3}:}
Let $I_{\F,\beta}(m,n)$ denote the maximum number of incidences between $m$ points and $n$ $\beta$-nondegenerate surfaces in $\F$. Fix some $\varepsilon>0$, we will prove by induction on $m+n$ that there exists some constant $C(\varepsilon, \beta, s, t)$ such that
\begin{equation}\label{eq:Ap-Sh for F}
I_{\F,\beta}(m,n)\leq C\left(m^{\frac{2s}{3s-1}+\varepsilon}n^{\frac{3s-3}{3s-1}}+mn^{1/2}+n\right)
\end{equation}
By choosing the constant big enough we can assume the bound holds for small values of $m+n$. For the induction step, let $j$ be a parameter to be chosen later. Let $\Gamma$ be the set of all intersection curves between surfaces in $\F$. A surface in $Q$ is called \emph{bad} if it contains a \emph{$j$-rich}  curve (i.e. a curve that contains at least $j$ points of $P$) from $\Gamma$,  and called \emph{good} otherwise. 
\\
\\
\textbf{Incidences from the good surfaces:}  For some parameter $r$ to be chosen later (think of $r$ as a big constant depending on $\varepsilon$), we construct a \emph{r-cutting} of the surfaces in $Q$ as follows: we choose each surface in $Q$ 
independently at random with probability $C_1 r/n\log r$ for an appropriate sufficiently large constant $C_1$. As shown in \cite{chazelle1991singly}, for any $\varepsilon_0>0$, we can construct a vertical decomposition of the arrangement of the random sample with $C_2r^{3+\varepsilon_0}$ open cells (where $C_2$ is some constant depending on $\varepsilon_0$), each of constant description complexity. With high probability, each cell is \emph{crossed}  (i.e., is intersected by but not contained in) by at most $\frac{n}{r}$ surfaces of $Q$. We may also assume that each cell contains at most $\frac{m}{r^3}$
points of $P$ by partitioning any cell which contains more into smaller cells (which can create at most $r^3$ new cells).

Consider a cell $\Omega$ and a good surface $S$ that crosses $\Omega$. If $S$ contains more than $\frac{1}{\beta}j$ points in $P\cap \Omega$ then $S$ is $\beta$-nondegenerate with respect to $P\cap \Omega$ (in fact it is geometrically $\beta$-nondegenerate). We can apply our induction assumption for the set of points in $\Omega$ and these spheres to get at most $I(\frac{m}{r^3}, \frac{n}{r})$ incidences. On the other hand, the contribution from the good surfaces with fewer than $\frac{1}{\beta}j$ points are $ \frac{jn}{\beta r}$. 
Therefore the total contribution from good surfaces is at most
$$C_2 r^{3+\varepsilon_0}\left(I\left(\frac{m}{r^3}, \frac{n}{r}\right)+\frac{jn}{\beta r}\right).$$

\noindent\textbf{Incidences from the bad surfaces:} Let $Q_0$ denote the set of all $j$-rich curves in $\Gamma$ which belong to some bad surface in $Q$. 
Project these curves to a generic plane, they still belong to a family of curves with $s-1$ degrees of freedom and multiplicity at most $t$ in a plane. By Pach-Sharir \cite{pach1998number}, the number of incidences between $P$ and $Q_0$ is at most $C_3\left(m^{\frac{s-1}{2(s-1)-1}}|Q_0|^{\frac{2(s-1)-2}{2(s-1)-1}}+m+|Q_0|\right)$ for some $C_3$ depending on $s,t$. On the other hand, since each curve is $j$-rich, there are at least $j|Q_0|$ incidences. We can assume $j>2C_3$ and get the following bound for some $C_4$ depending on $s,t$:
$$|Q_0|\leq C_4\left(\frac{m^{s-1}}{j^{2s-1}}+\frac{m}{j}\right).$$
For a fixed curve  $\gamma\in Q_0$, if it belongs to only one surface in $Q$, the total number of incidences involving $\gamma$ is at most $m$. Otherwise, when $\gamma$ belongs to at least two surfaces $S_1,S_2$, because of the nondegeneracy, $|\gamma\cap P|\leq \beta |S_i\cap P|$ for $i=1,2$.  This implies $|S\setminus\gamma\cap P|\geq (1-\beta)|S\cap P|$ for any surface $S$ containing $\gamma$. Moreover the sets in $\{S\setminus \gamma: S\in Q, \gamma\subset S\}$ are   mutually disjoint because any two such surfaces intersect at at most one irreducible curve, which is $\gamma$ in this case. Hence  
\begin{equation}\label{eq:disjoint argument}
\sum_{\gamma \subset S} |(S\setminus \gamma)\cap P|\leq m.
\end{equation}
So the total number of incidences in this case is at most 
$$|Q_0|\frac{1}{1-\beta}m\leq C_5 \left(\frac{m^s}{j^{2s-1}}+\frac{m^2}{j}\right)$$
for some $C_5$ depending on $s,t,\beta$.
Adding two cases we get
$$I_{\F, \beta}(m,n)\leq C_2 r^{3+\varepsilon_0}I_{\F,\beta}\left(\frac{m}{r^3},\frac{n}{r}\right)+C_2 r^{2+\varepsilon_0}jn+C_5\left(\frac{m^s}{j^{2s-1}}+\frac{m^2}{j}\right).$$
Applying the induction assumption for the first term and 
let $j=m^{\frac{1}{2}}n^{-\frac{1}{2s}}+mn^{-\frac{1}{2}}$ (to balance $jn$ with the last two terms), we get:
\begin{align}\label{eq:4 terms}
I_{\F, \beta}(m,n)&\leq C_2r^{3+\varepsilon_0} C\left(\left(\frac{m}{r^3}\right)^{\frac{2s}{3s-1}+\varepsilon}\left(\frac{n}{r}\right)^{\frac{3s-3}{3s-1}}+\frac{m}{r^3}\left(\frac{n}{r}\right)^{\frac{1}{2}}+n\right)+C_6r^{2+\varepsilon_0}\left(m^{\frac{1}{2}}n^{1-\frac{1}{2s}}+mn^{\frac{1}{2}}\right)\nonumber\\
&= CC_2r^{\varepsilon_0-3\varepsilon}m^{\frac{2s}{3s-1}+\varepsilon}n^{\frac{3s-3}{3s-1}}+(CC_2r^{\varepsilon_0-1/2}+C_6r^{2+\varepsilon_0})mn^{1/2}+Cr^{3+\varepsilon_0}n+C_6r^{2+\varepsilon_0}m^{\frac{1}{2}}n^{1-\frac{1}{2s}}.
\end{align}
We choose $\varepsilon_0<\min\{\varepsilon,1/4\}$ and $r$ big enough so that $\max\{C_2r^{\varepsilon_0-3\varepsilon},C_2r^{\varepsilon_0-1/2}\}<\frac{1}{2}$, but not too big so that $C_6r^{2+\varepsilon_0}<\frac{C}{4}$. This takes care of the first two terms. For the last two, we use \cref{thm:semi-alg graph} to get $I_{\F,\beta}=O_{\beta,s,t}(mn^{1-1/s}+n)$. Indeed, the incidence graph between $P$ and $Q$ is semi-algebraic with description complexity $t$ and $Q$ lives in a variety of dimension $s$\footnote{As observed in \cite{DO2018621}, it does not matter whether  $Q$ belongs to $\R^s$ or a $s$-dimensional variety of bounded degree in a larger space}. As a consequence, if $n>C_7m^{s}$ for some $C_7$, we have $I_{\F,\beta}=O_{\beta,s,t}(n)$, and we can choose $C$ large enough to make \eqref{eq:Ap-Sh for F} hold. So we can assume $n<C_7 m^s$. Simple algebra implies 
$$\max\{n,m^{\frac{1}{2}}n^{1-\frac{1}{2s}}\}=O\left( m^{\frac{2s}{3s-1}+\varepsilon}n^{\frac{3s-3}{3s-1}}\right).$$
This implies all four terms in \eqref{eq:4 terms} are bounded by the right hand side of \eqref{eq:Ap-Sh for F} and finishes the proof of \cref{thm:extend Ap-Sh in R3}.
\qed

\begin{proof}[Proof of \cref{cor:projected sphere}]
Project $P$ and $Q$ to a generic 3-dimensional flat $V$ to get $P',Q'$ respectively so that $P'$ is a set of $m$ distinct points,  $Q'$ consist of $n$ distinct ellipsoids,\footnote{The image of a sphere in this projection is either an ellipsoid or a solid ellipse. This is similar to if we project a circle in $\R^3$ to plane we either get an ellipse or a line segment.} and $I(P',Q')=I(P,Q)$. Let $\Gamma$ be the set of all intersection curves between two ellipsoids in $Q'$. Then it is easy to see each curve in $\Gamma$ is projection from some circle in $\R^4$, hence $\Gamma$ still   have three degrees of freedom. Each ellipsoid in $Q$ remains $\beta$-nondegenerate, i.e. contains no curve in $\Gamma$ with more than $\beta$ fraction of points. We would be done by applying \cref{thm:extend Ap-Sh in R3} if $Q'$ satisfies condition (ii) in \cref{def:good family}. Unfortunately this is not quite true: the parameter space for all ellipsoids has dimension six. However, if we look closely at the above proof, the only place we need (ii) is at the very end which allows us to use \cref{thm:semi-alg graph}. All we need is that $I(P',Q')=O_{\varepsilon}\left(m^{\frac{8}{11}+\varepsilon}n^{\frac{9}{11}}+mn^\frac{1}{2}\right)$ holds whenever $n>cm^4$ for some constant $c$.  

To show this, return to the space $\R^d$.  By \cref{lem:spanning}, $n\leq m^4$.
Consider the standard lifting: $\phi: \R^d\to\R^{d+1}$ where $\phi(x_1,\dots, x_d)=(x_1,\dots,x_d,\sum_{i=1}^d x_i^2)$. Each 2-sphere in $Q$ is mapped to a 3-flat in $\R^d$. Project   to a generic 4-flat $V$, these flats remain $\beta$-nondegenerate 3-flats in $H$. Hence we can apply Elekes-T\'oth's bound  \eqref{eq:E-T} in 4 dimension to get $I(P,Q)=O(m^{4/5}n^{4/5}+mn^{2/3})$. When $n>cm^4$ for some $c$, or $n=\Theta(m^4)$, this implies $I(P,Q)=O(n)=O(m^4)$ and the RHS of \eqref{eq:projected sphere} is $\Theta(m^4)$, hence \eqref{eq:projected sphere} holds true in this case. 
\end{proof}

\begin{proof}[Proof of \cref{lem:spanning}]
A hyperplane $H$ is called \emph{spanning} w.r.t. $P$ if $H$ contains some $d$ points in $P$ that are in general position (i.e. no three points collinear, no four points coplanar and so on). Equivalently, $H$ is spanning if it has no proper sub-flat that contains all the points of $H\cap P$. There are at most ${m\choose d}$ choices of $d$ points from $P$, each define at most one distinct spanning hyperplane, hence the number of spanning hyperplanes is at most ${m\choose d}$.
It is easy to see any geometrically $\beta$-nondegenerate hyperplane w.r.t. $P$ must be spanning, hence there are at most that many geometrically $\beta$-nondegenerate hyperplanes as well. 

The statement for spheres is proved similarly by using the concept of \emph{spanning spheres}, i.e.  one that contains $d+1$ points in $P$ that does not belong to any proper sub-sphere.
\end{proof}
\section{Point-spheres incidence graph in 4 dimensions}\label{sec:R4}
In this section we prove \cref{thm:nondeg spheres 4d}. Again we follow the strategy in \cite{apfelbaum2011non} with modifications. In particular (assuming the readers are familiar with the proof in the previous section), we need a good cutting method for spheres in $\R^4$ which is done by Koltun \cite{koltun2004almost}. The analysis for the bad 3-spheres, i.e. the ones that contain some $j$-rich 2-sphere, is more complicated: if the $j$-rich 2-spheres are $\beta$-nondegenerate, we can use \cref{cor:projected sphere} to bound the number of such spheres, then use the mutually disjoint argument in \eqref{eq:disjoint argument} . Otherwise they contain some $\beta j$-rich circles in which case we also get a bound on the number of such rich circles. What is significantly different in this case is that it is harder to bound the contribution from all the 3-spheres that contain a fixed circle. 
Here is the detailed proof.
\begin{proof}[Proof of \cref{thm:nondeg spheres 4d}] Let $I(m,n)$ denote the maximum number of incidences between $m$ points and $n$ 3-spheres in $\R^4$. 
We want to show that there exists a constant $C$ that depends on $\beta$ and $\varepsilon$ such that
\begin{equation}\label{eq:sphere 4d}
I(m,n)\leq C\left(m^{\frac{15}{19}+\varepsilon} n^{\frac{16}{19}}+mn^{\frac{2}{3}}\right).
\end{equation}

We prove by induction on $m+n$. By choosing the constant big enough we can assume the bound holds for small values of $m+n$. 
Let $j$ be a parameter to be chosen later. A 3-sphere is called \emph{bad} if it contains a $j$-rich two dimensional subsphere, and called \emph{good} otherwise. 
\\
\\
\textbf{Incidences from the good spheres:} We construct an \emph{r-cutting} of the spheres for some small $r$ to be chosen later. We sample spheres
by choosing each sphere in the sample
independently at random with probability $Cr/n\log r$ for an appropriate sufficiently large constant $C$. As shown in \cite{koltun2004almost}, for some small positive $\varepsilon_0$ to be chosen later, we can construct a vertical decomposition of the arrangement of the random sample with $O(r^{4+\varepsilon_0})$ open cells, each of constant description complexity. With high probability, each cell is crossed  (i.e., is intersected by but not contained in) by at most $\frac{n}{r}$ spheres of $Q$. We may also assume that each cell contains at most $\frac{m}{r^4}$
points of $P$ by partitioning any cell which contains more into smaller cells (which can create at most $r^4$ new cells).

Consider a cell $\Omega$ and a good 3-sphere $S$ that crosses $\Omega$. If $S$ contains more than $\frac{1}{\beta}j$ points in $P\cap \Omega$ then $S$ is $\beta$-nondegenerate with respect to $P\cap \Omega$. We can apply our induction assumption for the set of points in $\Omega$ and these spheres to get at most $I(\frac{m}{r^4}, \frac{n}{r}\log r)$ incidences. On the other hand, the contribution from the good spheres with fewer than $\frac{1}{\beta}j$ points are $O(j \frac{n}{r}\log r)$. Summing over all the cells, we get the total number of incidences with the good spheres is at most
\begin{equation}\label{eq:good sphere R4}
C_1 r^{4+\varepsilon_0}\left(I\left(\frac{m}{r^4}, \frac{n}{r}\right)+j\frac{n}{r}\right)
\end{equation}
for some constant $C_1$.
\\
\\
\noindent\textbf{Incidences from the bad spheres:} Let $Q_0$ denote the set of all $j$-rich 2-spheres which are subspheres of some bad spheres in $Q$. Divide $Q_0$ further into two sets: $Q_2$ which consist of all spheres that are $\beta$-nondegenerate with respect to $P$, and $Q_1$ the rest.  

By \cref{cor:projected sphere} we have
$|Q_2|\leq C_2\left(\frac{m^4}{j^{11/2}}+\frac{m^2}{j^2}\right)$ for some constant $C_2$.
For each sphere $\gamma\in Q_2$, the set $S\setminus \gamma$ for all spheres $S$ in $Q$ that contains $\gamma$ is disjoint. Hence  $\sum_{\gamma \subset S} |S\setminus \gamma|\leq m$. Because of the nondegenerate condition $|S\setminus\gamma|\geq (1-\beta)|S|$. The total number of incidences from spheres in $Q$ that contain some subsphere in $Q_2$ are then
\begin{equation}\label{eq:bad sphere Q2}
\sum_{\gamma\in Q_2}\sum_{\S\supset \gamma}|S|\leq |Q_2|\frac{1}{1-\beta}m\leq \frac{C_2}{1-\beta}\left(\frac{m^5}{j^{11/2}}+\frac{m^3}{j^2}\right).
\end{equation}
Each 2-sphere in $Q_1$ contains some  $(\beta j)$-rich circle. By Pach-Sharir's bound in \cite{pach1998number}, the number of $(\beta j)$-rich circles is\footnote{There is a stronger point-circle incidence bound in \cite{marcus2006intersection}, but this one has fewer terms and is enough for our purpose.} 
$$O\left(\frac{m^3}{(\beta j)^{5}}+\frac{m}{\beta j}\right).$$
For each such circle $C$, let $\Gamma_C$ denote the set of all 2-spheres in $Q_2$ that degenerates to $C$.
\begin{align*}
\sum_{C\subset S\in Q}|S\setminus C| & = \sum_{\gamma\in \Gamma_C}\sum_{S\supset\gamma}|S\setminus C|\\
&=\sum_{\gamma\in \Gamma_C}\sum_{S\supset\gamma}\left(|S\setminus \gamma|+|\gamma\setminus C|\right)\\
&=\sum_{\gamma\in \Gamma_C}\sum_{S\supset\gamma}|S\setminus \gamma|+\sum_{S\supset\gamma}\sum_{\gamma\in \Gamma_C}|\gamma\setminus C|\\
&\leq m+\sum_{S\supset\gamma} m\\
&\leq m+\frac{m}{(1-\beta)j}m
\end{align*}
In the second the last inequality, we use the fact that all the sets in each summation are disjoint (similarly to the argument for $Q_2$ above). In the last inequality, we use the fact that each $\gamma\in Q_1$ is belonged to at most $\frac{m}{(1-\beta)j}$ 3-spheres in $Q$ because $S\setminus \gamma$ are all disjoint and each contains at least $(1-\beta)j$ points.

Since $|C|\geq \beta|\gamma|\geq \beta^2|S|$, we have again $|S|\leq \frac{1}{1-\beta^2}|S\setminus C|$ for each sphere $S\in Q$ that contains $C$. So the total number of incidences from spheres in $Q$ that contain some 2-spheres in $Q_1$ is at most
\begin{equation}\label{eq:bad sphere Q1}
C_3\left[\left(\frac{m^2}{j}+m\right)\left(\frac{m^3}{(\beta j)^{5}}+\frac{m}{\beta j}\right)\right]=C_3\left(\frac{m^5}{j^{6}}+\frac{m^4}{j^5}+\frac{m^3}{j^2}+\frac{m^2}{j}\right)
\end{equation}
for some constant $C_3$. Combining \eqref{eq:good sphere R4},\eqref{eq:bad sphere Q2} and \eqref{eq:bad sphere Q1} we have
$$I(m,n)\leq C_1 r^{4+\varepsilon_0}I\left(\frac{m}{r^4}, \frac{n}{r}\right)+C_1jnr^{3+\varepsilon_0}+ \frac{C_2}{1-\beta}\left(\frac{m^5}{j^{11/2}}+\frac{m^3}{j^2}\right)+C_3\left(\frac{m^5}{j^{6}}+\frac{m^4}{j^5}+\frac{m^3}{j^2}+\frac{m^2}{j}\right)$$
$$\leq C_1 r^{4+\varepsilon_0}I\left(\frac{m}{r^4}, \frac{n}{r}\right)+C_4\left(r^{3+\varepsilon_0} nj+\frac{m^5}{j^{11/2}}+\frac{m^4}{j^5}+\frac{m^3}{j^2}+\frac{m^2}{j}\right).$$
Apply the induction assumption for the first term: $I(\frac{m}{r^4},\frac{n}{r})\leq C\left(\frac{m}{r^4}\right)^{\frac{15}{19}+\varepsilon}\left(\frac{n}{r}\right)^\frac{16}{19}+C\frac{m}{r^4}(\frac{n}{r})^{\frac{2}{3}}$ and let $j=m^{10/13}n^{-2/13}+m^{4/6}n^{-1/6}+mn^{-1/3}+mn^{-1/2}$ to balance $nj$ with the rest, we get
$$I(m,n)=O_{\varepsilon, \varepsilon_0}\left[r^{\varepsilon_0-4\varepsilon}m^{\frac{15}{19}+\varepsilon}n^{\frac{16}{19}}+r^{\varepsilon_0-2/3} mn^{\frac{2}{3}}+r^{3+\varepsilon_0}\left(m^\frac{10}{13}n^\frac{11}{13}+m^\frac{4}{6}n^\frac{5}{6}+mn^\frac{2}{3}+mn^\frac{1}{2} \right)\right].$$
Similarly to the proof in the previous section, we can choose $\varepsilon_0$ so that the exponents of $r$ in the first two terms are negative, hence we can choose $r$ so that the constants in the first two term is less than $1/2 C$. It is easy to see $m^\frac{4}{6}n^\frac{5}{6}<m^{\frac{15}{19}}n^{\frac{16}{19}}$ and $mn^{1/2}<mn^{2/3}$. Finally, we have $m^\frac{10}{13}n^\frac{11}{13}<m^{\frac{15}{19}}n^{\frac{16}{19}}$ when $n<m^5$. We can assume $n\leq m^5$ because of \cref{lem:spanning}.
So if we choose $r$ such that $r^{3+\varepsilon_0}<m^\varepsilon$, we can bound the all the terms by the RHS of \eqref{eq:sphere 4d}. This finishes the proof of \cref{thm:nondeg spheres 4d}.
\end{proof}

\section{Application to the similar simplices problem}\label{section:application}
In this section we will prove \cref{thm:application to similar triangles in R4}. We follow the strategy in \cite{agarwal2007similar} for bounding the number of similar triangles in three dimensions, but again the details are much more complicated with one extra dimension.
\begin{proof}[Proof of \cref{thm:application to similar triangles in R4}]
Fix two points $a,b$ in $P$. The orbit of all points $c$ that form a triangle similar to $\Delta$ where $ab$ is the longest edge and $ac$ is the second longest one\footnote{If $\Delta$ has more than one longest edges, pick either one and the calculation may be off by at most a constant factor. } is the intersection of two 3-spheres centered at $a,b$, which is a 2-sphere lying on a hyperplane perpendicular to the line spanned by $a$ and $b$.  We denote this 2-sphere by $\gamma_{a,b}$, and let $\Gamma:=\{\gamma_{a,b}:a,b\in P\}$. The number of triangles similar to $\Delta$ is of the same order with the number of incidences between $n$ points in $P$ and $O(n^2)$ 2-spheres in $\Gamma$. The result in \cite{agarwal2007similar} for $f_{2,4}(n)$ is obtained by estimating this incidence which depends on whether these 2-spheres are nondegenerate or not.
We, on the other hand, use the fact that the point $c$ also belongs to a 3-sphere $S$ centered at $a$, hence reduce the problem to estimate incidences between points and 2-spheres $\gamma_{a,b}$ within a 3-sphere $S$. This is similar to what was used to prove the bound for the number of similar triangles in $\R^3$ in \cite{apfelbaum2011non}.

For a fixed vertex $a\in P$, let $R_a$ denote all possible distances from other points in $P$ to $a$. For each $r\in R$, let $S_{a,r}$ be  the 3-sphere centered at $a$ and radius $r$ and $|S_{a,r}|$ be the number of points of $P$ which lie in $S_{a,r}$. Clearly $\sum_{r\in R_a}|S_{a,r}|=n-1$. 

The number of triangles similar to $\Delta$ with $a$ as a vertex adjacent to two longest edge and $b\in S_{a,r}$ is equal to the number of incidences between $|S_{a,r}|$ 2-spheres $\{\gamma_{a,b}:b\in S_{a,r}\}$ (of same radius) and $|S_{a,r\cdot d}|$ points in $S_{a,r\cdot d}$ where $d$ is the ratio between the second longest edge with the longest one in $\Delta$.  By Theorem 6.4 in \cite{fox2014semi}, the number of incidences between $n$ unit spheres and $m$ points on a 3-sphere is $O((mn)^{3/4+\varepsilon}+m+n)$ as their incidence graph is $K_{3,3}$-free (which holds true in $\R^3$ and also in a 3-sphere). Therefore we have the following bound:
$$I\lesssim \sum_{a\in P}\sum_{r\in R_a} \left((|S_{a,r}||S_{a,r\cdot d}|)^{3/4+\varepsilon}+|S_{a,r}|+|S_{a,r\cdot d}|\right)\lesssim \sum_{a\in P}\sum_{r\in R_a} (|S_{a,r}||S_{a,r\cdot d}|)^{3/4+\varepsilon}+ 2n(n-1).$$
We only need to bound the first term. Let us consider three cases: (1) either $|S_{a,r}|$ or $|S_{a,r\cdot d}|$ is at most $t$ for some parameter $t$ to be chosen later, (2) both $S_{a,r}$ and $S_{a,dr}$ are $t$-rich and  $\beta$-nondegenerate for some $0<\beta<1$, and (3) both $S_{a,r}$ and $S_{a,dr}$ are $t$-rich and one of them is $\beta$-degenerate.   Let $I_1, I_2, I_3$ denote the contribution to $I$ from case 1,2,3 respectively, then $I=I_1+I_2+I_3$.
\\
\\
\textbf{Bounding $I_1$:} 
WLOG we can assume $|S_{a,r}|\leq t$. Using H\"older's inequality the contribution from this case is 
\begin{align*}I_1\leq & \sum_{a\in R} t^{1/2+2\varepsilon}\sum_{r\in R} |S_{a,r}|^{1/4-\varepsilon}|S_{a,r\cdot d}|^{3/4+\varepsilon}\\
\leq & \sum_{a\in R} t^{1/2+2\varepsilon} \left(\sum_{r\in R_a} |S_{a,r}|\right)^{1/4+\varepsilon}\left(\sum_{r\in R_a} |S_{a,r\cdot d}|\right)^{3/4+\varepsilon}\\
\leq & \sum_{a\in P} t^{1/2+\varepsilon} (n-1)\\
\leq & t^{1/2+\varepsilon} n^2
\end{align*}

\noindent\textbf{Bounding $I_2$:} 
We use dyadic decomposition: let $R_i$ denote the set of all pairs $(a,r)$ where $r\in R_a$ such that 
$$2^{i}t<\max\{|S_{a,r}|,|S_{a,dr}|\}<2^{i+1} t.$$
By \cref{thm:nondeg spheres 4d}, the number of $2^{i}t$-rich $\beta$-nondegenerate 3-spheres with respect to $P$ is $O(\frac{n^{5+\varepsilon}}{(2^i t)^{19/3}}+\frac{n^3}{(t2^i)^3})$.
\begin{align*}
I_2\lesssim & \sum_{i=1}^{\log (n/t)}\sum_{(a,r)\in R_i} (|S_{a,r}||S_{a,dr}|)^{3/4+\varepsilon}\\
\leq & \sum_{i=1}^{\log (n/t)} |R_i| (2^{i+1}t)^{3/2+2\varepsilon}\\
\lesssim & \sum_{i=1}^{\log (n/t)} \left(\frac{n^4}{(2^i t)^{19/3}}+\frac{n^3}{(t2^i)^3}\right) (2^{i+1}t)^{3/2+2\varepsilon}\\
\leq & \frac{n^{5+\varepsilon}}{t^{29/6-2\varepsilon}}\sum_{i=1}^{\log (n/t)}\frac{1}{2^{i(29/6-2\varepsilon)}}+\frac{n^3}{t^{3/2-2\varepsilon}}\sum_{i=1}^{\log (n/t)}\frac{1}{2^{3/2 i}}\\
\lesssim & \frac{n^{5+\varepsilon}}{t^{29/6-2\varepsilon}}+\frac{n^3}{t^{3/2-2\varepsilon}}
\end{align*}

\noindent\textbf{Bounding $I_3$:} Both $S_{a,r}$ and $S_{a,dr}$ are $t$-rich and one of them is $\beta$-degenerate. WLOG assume $S_{a,dr}$ is $\beta$-degenerate.\footnote{The other case can be treated similarly by rewriting $I_3$ as incidences between $|S_{a,dr}|$ 2-spheres and $|S_{a,r}|$ points  on $S_{a,r}$.} This means there exists a 2-sphere $\gamma_1$ inside $S_{a,dr}$ which contains at least $\beta|S_{a,dr}|$ points of $P$. We remove the points of $P\cap \gamma_1$. If $S_{a,dr}$ remains $t$-rich $\beta$-degenerate, find another rich 2-sphere $\gamma_2$ and remove it. Keep doing it until $S_{a,dr}$ is either non $t$-rich then apply case 1, or $\beta$-nondegenerate then apply case 2. Hence it remains to estimate the contribution to $I$ from the removed points. 

There are two cases. First consider incidences of the form $(p,\gamma)$ where $p$ is removed and $\gamma$ is a removed 2-sphere. For any  $\gamma$, there are at most two choices for $(a,b)$ such that $\gamma=\gamma_{a,b}$, hence $\gamma$ can be removed at most twice. The contribution from this type of incidence is thus at most twice the number of incidences between $P$ and the set of $\beta t$-rich spheres in $\Gamma$. If the spheres are $\beta$-nondegenerate, by \eqref{eq:Ap-Sh} and dyadic decomposition, the number of incidences is at most
$$O\left(\frac{n^4}{t^{9/2}}+\frac{n^2}{t}\right).$$
If the spheres are $\beta$-degenerate, each contains a $\beta^2 t$-rich circles. The number of point-circle incidences with these rich circles is at most $O(\frac{n^3}{t^{9/2}}+\frac{n^2}{t^2}+n)$ by bounding the number of rich circles and dyadic decomposition. Moreover, each circle has multiplicity at most $n$ because for each choice of $a\in P$, there is at most one choice for $b$ so that $\gamma_{a,b}$ contains this circle: the point on a line through $a$ and the center of the circle and has a certain distance to $a$.  So the total contribution in this case is 
$O\left(\frac{n^4}{t^{9/2}}+\frac{n^3}{t^2}+n^2\right).$

Now consider incidences of form $(p, \gamma)$ where $p$ is removed and $\gamma$ is not a removed 2-sphere. Note that the total number of removed 2-spheres on $S_{a,dr}$ is $O(\log n)$ because each time removing a sphere, the number of points shrinks by at least $1-\beta$ fraction. Hence for a fixed non-removed sphere $\gamma$, any point $p$ being removed must belong to $O(\log n)$ circles which are the intersection of $\gamma$ with the removed spheres. The contribution to $I_3$ from such circles that are not $\sqrt{t}$-rich is $O(\sqrt{t} n^2\log n )$, while the contribution from  $\sqrt{t}$-rich circles, as argued in the previous paragraph, is $$O\left(\frac{n^4}{(\sqrt{t})^{9/2}}+\frac{n^3}{(\sqrt{t})^2}+n^2\right).$$

Summing all three cases we have
$$I\leq  t^{1/2+\varepsilon} n^2+\frac{n^{5+\varepsilon}}{t^{29/6-2\varepsilon}}+\frac{n^3}{t^{3/2-2\varepsilon}}+n^2\log n \sqrt{t}+\frac{n^4}{(\sqrt{t})^{9/2}}+\frac{n^3}{(\sqrt{t})^2}+n^2.$$
Finally, by choosing $t=n^{8/11}$, we get the bound $I=O(n^{2+4/11+\varepsilon})$.
\end{proof}
\bibliographystyle{abbrv}
\bibliography{non-deg-sphere.bib}

\end{document}